\input ./style/arxiv-general.cfg
\documentclass[MSNbibl,number,citesort,seceqn,dvips]{arxbj}
\makeatletter
   \@ifpackageloaded{graphicx}{}{\usepackage{graphicx}}
\makeatother

\usepackage{mathrsfs}

\volume{22}
\issue{4}
\pubyear{2016}
\firstpage{1979}
\lastpage{2000}
\doi{10.3150/15-BEJ717}
\docsubty{FLA}

\makeatletter
\newcommand{\rrvert}{\vert}
\newcommand{\rrVert}{\Vert}
\newcommand{\llvert}{\vert}
\newcommand{\llVert}{\Vert}
\newtheorem{theorem}{Theorem}[section]
\newtheorem{proposition}[theorem]{Proposition}
\newtheorem{corollary}[theorem]{Corollary}
\newremark{remark}[theorem]{Remark}
\newcommand{\BX}{\mathbf{X}}
\newcommand{\reals}{{\mathbb R}}
\newcommand{\bbn}{{\mathbb N}}
\newcommand{\bbs}{{\mathbb S}}
\newcommand{\sG}{{\mathcal G}}
\newcommand{\one}{\mathbf{1}}
\makeatother

\begin{document}
\begin{frontmatter}

\title{Time-changed extremal process as a random sup measure}
\runtitle{Time-changed extremal process}

\begin{aug}
\author[A,B,C]{\inits{C.}\fnms{C\'eline}~\snm{Lacaux}\corref{}\thanksref{A,B,C}\corref{}\ead[label=e1]{Celine.Lacaux@univ-lorraine.fr}}
\and
\author[D]{\inits{G.}\fnms{Gennady}~\snm{Samorodnitsky}\thanksref{D}\ead[label=e2]{gs18@cornell.edu}}
\address[A]{Universit\'e de Lorraine, Institut \'Elie Cartan de
Lorraine, UMR 7502, Vand\oe uvre-l\`es-Nancy, F-54506, France.
\printead{e1}}
\address[B]{CNRS, Institut \'Elie Cartan de Lorraine, UMR 7502,
Vand\oe uvre-l\`es-Nancy, F-54506, France}
\address[C]{Inria, BIGS, Villers-l\`es-Nancy, F-54600, France}
\address[D]{School of Operations Research and Information Engineering
and Department of Statistical Science Cornell University, Ithaca, NY
14853, USA. \printead{e2}}
\end{aug}

%
\received{\smonth{10} \syear{2014}}
%
\revised{\smonth{2} \syear{2015}}

%
\begin{abstract}
A functional limit theorem for the partial maxima of a long memory
stable sequence produces a limiting process that can be described as a
$\beta$-power time change in the classical Fr\'echet
extremal process, for $\beta$ in a subinterval of the unit
interval. Any such power time change in the extremal process
for $0<\beta<1$ produces a process with stationary
max-increments. This deceptively simple time change hides the much
more delicate structure of the resulting process as a self-affine
random sup measure. We uncover this structure and show that in a
certain range of the parameters this random measure arises as a limit
of the partial maxima of the same long memory stable sequence, but in
a different space. These results open a way to construct a whole new
class of self-similar Fr\'echet processes with stationary
max-increments.
\end{abstract}

%
\end{frontmatter}

\section{Introduction} \label{secintro}

Let $(X_1,X_2,\ldots)$ be a stationary sequence of random variables,
and let $M_n=\max_{1 \leq k \leq n} X_k$, $n=1,2,\dots$ be the
sequence of its partial maxima. The limiting distributional behavior
of the latter sequence is one of the major topics of interest in extreme
value theory. We are particularly interested in the possible limits in a
functional limit theorem of the form
%
%
\begin{equation}
\label{efunctconv} \biggl( \frac{M_{\lfloor nt\rfloor}-b_n}{a_n}, t\geq0 \biggr) \Rightarrow\bigl(Y(t), t
\geq0\bigr),
\end{equation}
for properly chosen sequences $(a_n)$, $(b_n)$. The weak convergence in
{(\ref{efunctconv})} is typically in the space $D[0,\infty)$ with one
of the usual Skorohod topologies on that space; see
\cite{skorohod1956,billingsley1999} and
\cite{whitt2002}. If the original sequence $(X_1,X_2,\ldots)$ is an
i.i.d. sequence, then the only possible limit in~(\ref{efunctconv})
is \textit{the extremal process}, the extreme value analog of the
L\'evy process; see \cite{lamperti1964}.

The modern extreme value theory is interested in
the case when the sequence $(X_1,X_2,\ldots)$ is
stationary, but not necessarily independent. The potential clustering
of the extremes in this case leads one to expect that new limits may
arise in~(\ref{efunctconv}). Such new limits, however, have not
been widely observed, and the dependence in the model has been
typically found to be reflected in the limit via a linear time change
(a slowdown), often connected to the \textit{extremal
index}, introduced originally in \cite{leadbetter1983}. See,
for example, \cite{leadbetterlindgrenrootzen1983}, as well as the studies
in \cite{rootzen1978,davisresnick1985,mikoschstarica2000b} and \cite{fasen2005}.
One possible explanation for this is the known phenomenon that the
operation of taking partial maxima tends to mitigate the effect of
dependence in the original stationary sequence, and the dependent
models considered above were, in a certain sense, not sufficiently
strongly dependent.

Starting with a long range dependent sequence may make a difference,
as was demonstrated by \cite{owadasamorodnitsky2014}. In that paper,
the original sequence was (the absolute value of) a stationary
symmetric $\alpha$-stable process, $0<\alpha<2$, and the length of
memory was quantified by a single parameter $0<\beta<1$. In the case
$1/2<\beta<1$, it was shown that the limiting process in
(\ref{efunctconv}) can be represented in the form
%
%
\begin{equation}
\label{elimprocess} Z_{\alpha, \beta}(t) = Z_{\alpha}\bigl(t^\beta
\bigr),\qquad t\geq0,
\end{equation}
where $ ( Z_{\alpha}(t),   t\geq0 )$ is the extremal
($\alpha$-)Fr\'echet process.

The nonlinear power time change in~(\ref{elimprocess}) is both
surprising and misleadingly simple. It is surprising because it is not
immediately clear that such a change is compatible with a certain
translation invariance the limiting process must have due to the
stationarity of the original sequence. It is misleadingly simple
because it hides a much more delicate structure. The main goal of this
paper is to reveal that structure. We start by explaining exactly what
we are looking for.

The stochastic processes in the left-hand side of
(\ref{efunctconv}) can be easily interpreted as \textit{random sup
measures} evaluated on a particular family of sets (those of the form
$[0,t]$ for $t\geq0$). If one does not restrict himself to that specific
family of sets and, instead, looks at all Borel subsets of
$[0,\infty)$, then it is possible to ask whether there is weak
convergence in the appropriately defined space of random sup
measures, and what might be the limiting random sup measures. See the
discussion around~(\ref{epartialmn}) and the convergence result in
Theorem~\ref{textendOS}. This is
the approach taken in~\cite{obrientorfsvervaat1990}. Completing the
work published in \cite{vervaat1986} and \cite{vervaat1997}, the
authors provide a detailed description of the possible limits. They
show that the limiting random sup measure must be \textit{self-affine}
(they refer to random sup measures as extremal processes, but we
reserve this name for a different object).

As we will see in the sequel, if~(\ref{efunctconv}) can be stated in
terms of weak convergence of a sequence of random
sup measures, this would imply the finite-dimensional convergence part
in the functional formulation of~(\ref{efunctconv}). Therefore, any
limiting process $Y$ that can be obtained as a limit in this case must
be equal in
distribution to the restriction of a random sup measure to the sets of
the form
$[0,t]$, $t\geq0$. The convergence to the process $Z_{\alpha, \beta}$
established in \cite{owadasamorodnitsky2014} was not established in
the sense of weak convergence of a sequence of random
sup measures, and one of our tasks in this paper is to fill this gap and
prove the above convergence. Recall, however, that the convergence in
\cite{owadasamorodnitsky2014} was established only for $0<\alpha<2$
(by necessity, since $\alpha$-stable processes do not exist outside of
this range) and $1/2<\beta<1$. The nonlinear time change in
(\ref{elimprocess}) is, however, well defined for all $\alpha>0$ and
$0<\beta<1$, and leads to a process $Z_{\alpha, \beta}$ that is
self-similar and has stationary max-increments. Our second task in
this paper is to prove that the process $Z_{\alpha, \beta}$ can, for
\textit{all values of its parameters}, be extended to a random sup measure
and elucidate the structure of the resulting random sup measure. The
key result is Corollary~\ref{cmainresult} below.
The structure we obtain is of interest on its own right. It is
constructed based
on a certain random closed set possessing appropriate scaling
and translation invariance properties. Extending this approach
to other random sets and other ways of handling these random sets
may potentially lead to a construction of new classes
of self-similar processes with stationary max-increments and of random
sup measures. This is important both theoretically, and may be useful
in applications.

This paper is organized as follows. In the next section, we will define
precisely the notions discussed somewhat informally above and
introduce the required technical background. Section~\ref{secergodic}
contains a discussion of the dynamics of the stationary sequence
considered in this paper. It is based on a null recurrent Markov chain.
In Section~\ref{secsupmeasure}, we will prove that the process
$Z_{\alpha,
\beta}$ can be extended to a random sup measure and construct
explicitly such an extension. In Section~\ref{secconvergence}, we
show that the convergence result of \cite{owadasamorodnitsky2014}
holds, in a special case of a Markovian ergodic system, also in the
space SM of sup measures. Finally, in Section~\ref{secotherpr} we
present one of the possible extensions of the present work.

\section{Background} \label{secbackground}

An extremal process $ ( Y(t),  t\geq0 )$ can be viewed as an
analog of a L\'evy motion when the operation of summation is replaced
by the operation of taking the maximum. The one-dimensional marginal
distribution of a L\'evy process at time 1 can be an arbitrary infinitely
divisible distribution on $\reals$; \textit{any} one-dimensional
distribution is infinitely divisible with respect to the operation of
taking the maximum. Hence the one-dimensional marginal
distribution of an extremal process at time 1 can be any distribution
on $[0,\infty)$; the restriction to the nonnegative half-line being
necessitated by the fact that, by convention, an extremal process,
analogously to
a L\'evy process, starts at the origin at time zero. If $F$ is the
c.d.f. of a
probability distribution on $[0,\infty)$, then the finite-dimensional
distributions of an extremal process with distribution $F$ at time 1
can be defined by
%
%
\begin{eqnarray}\label{eextremepr}
&& \bigl( Y(t_1),Y(t_2),\ldots,
Y(t_n) \bigr)
\nonumber\\[-8pt]\\[-8pt]\nonumber
&&\quad \stackrel{d} {=} \bigl( X^{(1)}_{t_1},
\max \bigl( X^{(1)}_{t_1}, X^{(2)}_{t_2-t_1}
\bigr), \ldots,
\max \bigl( X^{(1)}_{t_1}, X^{(2)}_{t_2-t_1},
\ldots, X^{(n)}_{t_n-t_{n-1}} \bigr) \bigr)
\nonumber
\end{eqnarray}
for all $n\ge1$ and $0\leq t_1<t_2<\cdots<t_n$. The different random
variables in the right-hand side of~(\ref{eextremepr}) are
independent, with $X^{(k)}_{t}$ having the c.d.f. $F^t$ for $t>0$. In
this paper, we deal with the $\alpha$-Fr\'echet extremal process, for
which
%
%
\begin{equation}
\label{efrechet} F(x) = F_{\alpha,\sigma}(x) = \exp \bigl\{ -\sigma^\alpha
x^{-\alpha} \bigr\},\qquad x>0,
\end{equation}
the Fr\'echet law with the tail index $\alpha>0$ and the scale
$\sigma>0$. A stochastic process $(Y(t),   t \in T)$ (on an
arbitrary parameter
space $T$) is called a Fr\'echet process if for all $n \geq1$, $a_1,
\dots, a_n >0$ and $t_1, \dots, t_n \in T$, the weighted maximum
$\max_{1 \leq j \leq n} a_j Y(t_j)$ has a Fr\'echet law as in
(\ref{efrechet}). Obviously, the Fr\'echet extremal process is an
example of a Fr\'echet process, but there are many Fr\'echet
processes on $[0,\infty)$ different from the Fr\'echet extremal
process; the process $Z_{\alpha, \beta}$ in~(\ref{elimprocess}) is
one such process.

A stochastic process $ ( Y(t),
  t\geq0 )$ is called self-similar with exponent $H$ of
self-similarity if for any $c>0$
\[
\bigl( Y(ct), t\geq0 \bigr)\stackrel{d} {=} \bigl( c^HY(t), t\geq0
\bigr) %
\]
in the sense of equality of finite-dimensional distributions.
A stochastic process $( Y(t),   t \geq0 )$ is said to have
stationary max-increments if for every $r \geq0$, there exists,
perhaps on an enlarged probability space, a
stochastic process $ ( Y^{(r)}(t),   t \geq0  )$ such that
%
\begin{equation}
\label{edefstatmaxi} \cases{ \bigl( Y^{(r)}(t), t \geq0 \bigr) \stackrel{d}
{=} \bigl( Y(t), t \geq0 \bigr), 
\vspace*{3pt}\cr
\bigl( Y(t+r), t \geq0 \bigr)
\stackrel{d} {=} \bigl( Y(r) \vee Y^{(r)}(t), t \geq0 \bigr),}
\end{equation}
%
with $a\vee b=\max(a,b)$; see \cite{owadasamorodnitsky2014}.
This notion is an analog of the usual notion of a process with stationary
increments (see, e.g., \cite{embrechtsmaejima2002} and
\cite{samorodnitsky2006LRD})
suitable for the situation where the operation of summation
is replaced by the operation of taking the maximum. It follows from
Theorem 3.2 in \cite{owadasamorodnitsky2014} that only self-similar
processes with stationary max-increments can be obtained as limits in
the functional convergence scheme~(\ref{efunctconv}) with $b_n\equiv
0$.

We switch next to a short overview of random sup measures. The reader
is referred to \cite{obrientorfsvervaat1990} for full details. Let
$\sG$ be the collection of open subsets of $[0,\infty)$. We call a map
$m:  \sG\to[0,\infty]$ a sup measure (on $[0,\infty)$) if
$m(\varnothing)=0$ and
\[
m \biggl( \bigcup_{r\in R} G_r \biggr) =
\sup_{r\in R} m(G_r) %
\]
for an arbitrary collection $ ( G_r,   r\in R )$ of open
sets. In general, a sup measure can take values in any closed
subinterval of $[-\infty,\infty]$, not necessarily in $[0,\infty]$,
but we will consider, for simplicity, only the nonnegative case in
the sequel, and restrict ourselves to the maxima of nonnegative random
variables as well.

The \textit{sup derivative} of a sup measure is a function
$[0,\infty)\to[0,\infty]$ defined by
\[
d \, \check{\phantom{_i}} m(t) = \inf_{ G\ni t} m(G),\qquad t\geq0. %
\]
It is automatically an upper semicontinuous function. Conversely, for
any function $f:  [0,\infty)\to[0,\infty]$ the \textit{sup integral} of
$f$ is a sup measure defined by
\[
i \, \check{\phantom{_i}} f(G) = \sup_{t\in G}f(t),\qquad G\in\sG, %
\]
with $i \, \check{\phantom{_i}} f(\varnothing)=0$ by convention.
It is always true that $m=i \, \check{\phantom{_i}} d \, \check{\phantom{_i}} m$ for any sup
measure $m$, but
the statement $f=d \, \check{\phantom{_i}} i \, \check{\phantom{_i}} f$ is true only for upper
semicontinuous
functions $f$. A sup measure has a canonical extension to all subsets of
$[0,\infty)$ via
\[
m(B) = \sup_{t\in B} d \, \check{\phantom{_i}} m(t). %
\]
On the space SM of sup measures, one can introduce a topology, called
the \textit{sup vague topology} that makes SM a compact metric space. In
this topology, a sequence $(m_n)$ of sup measures converges to a sup
measure $m$ if both
\[
\limsup_{n\to\infty} m_n(K) \leq m(K)\qquad\mbox{for
every compact }K %
\]
and
\[
\liminf_{n\to\infty} m_n(G) \geq m(G)\qquad\mbox{for
every open }G. %
\]
A random sup measure is a measurable map from a probability space into
the space SM equipped with the Borel $\sigma$-field generated by the sup
vague topology.

The convergence scheme~(\ref{efunctconv}) has a natural version in
terms of random sup measures. Starting with a stationary sequence
$\BX=(X_1,X_2,\ldots)$ of nonnegative random variables, one can
define for
any set $B\subseteq[0,\infty)$
%
%
\begin{equation}
\label{epartialmn} M_n(\BX) (B) = \max_{k:  k/n\in B}
X_k.
\end{equation}
Then for any $a_n>0$, $M_n(\BX)/a_n$ is a random sup measure, and
\cite{obrientorfsvervaat1990} characterizes all possible limiting
random sup measures in a statement of the form
%
%
\begin{equation}
\label{esupmconv} \frac{M_n(\BX)}{a_n} \Rightarrow M
\end{equation}
for some sequence $(a_n)$. The convergence is weak convergence in the
space SM equipped with the sup vague topology. Theorem 6.1 in \cite{obrientorfsvervaat1990}
shows that any limiting random sup measure $M$ must be both
stationary and self-similar, that is,
%
%
\begin{equation}
\label{esasupm} M(a+\cdot)\stackrel{d} {=}M\quad\mbox{and}\quad
a^{-H}M(a\cdot) \stackrel{d} {=}M\qquad\mbox{for all }a>0
\end{equation}
for\vspace*{1pt} some exponent $H$ of self-similarity. In fact, the results of
\cite{obrientorfsvervaat1990} allow for a shift $(b_n)$ as in
(\ref{efunctconv}), in which case the power scaling $a^{-H}$ in
(\ref{esasupm}) is, generally, replaced by the scaling of the form
$\delta^{-\log a}$, where $\delta$ is an affine transformation. In the
context of the present paper, this additional generality does not play
a role.

Starting with a stationary and self-similar random sup measure $M$,
one defines a stochastic process by
%
%
\begin{equation}
\label{erestrM} Y(t) = M \bigl( (0,t] \bigr),\qquad t\geq0.
\end{equation}
Then the self-similarity property of the random sup measure $M$
immediately implies the self-similarity property of the stochastic
process $Y$, with the same exponent of self-similarity. Furthermore,
the stationarity of the random sup measure $M$ implies that the
stochastic process $Y$ has stationary max-increments; indeed, for
$r\geq0$ one can simply take
\[
Y^{(r)}(t) = M \bigl( (r,r+t] \bigr),\qquad t\geq0. %
\]

Whether or not any self-similar process with stationary
max-increments can be constructed in this way or, in other words,
whether or not such a process can be extended, perhaps on an extended
probability space, to a stationary and
self-similar random sup measure remains, to the best of our knowledge,
an open question. We do show that the process $Z_{\alpha, \beta}$ in
(\ref{elimprocess}) has such an extension.

\section{The Markov chain dynamics}
\label{secergodic}

The stationary sequence we will consider in Section~\ref{secconvergence} is a symmetric $\alpha$-stable (S$\alpha$S) sequence,
whose dynamics is driven by a certain Markov chain. Specifically,
consider an irreducible null recurrent Markov chain $(Y_n,   n \geq
0)$ defined on an infinite countable state space $\bbs$ with
transition matrix $(p_{ij})$. Fix an arbitrary state $i_0\in\bbs$,
and let $(\pi_i,   i \in\bbs)$ be the unique invariant measure of
the Markov chain with
$\pi_{i_0}=1$. Note that $(\pi_i)$ is necessarily an infinite
measure.

Define a $\sigma$-finite and infinite measure on
$(E,\mathcal{E}) = (\bbs^{\bbn}, \mathcal{B}(\bbs^{\bbn}))$ by
\[
\mu(B) = \sum_{i \in\bbs} \pi_i
P_i(B),\qquad B \in\mathcal{E}, %
\]
where $P_i(\cdot)$ denotes the probability law of $(Y_n)$ starting in
state $i \in\bbs$. Clearly, the usual left shift operator on
$\bbs^{\bbn}$
\[
T(x_0, x_1, \dots) = (x_1,x_2,
\dots) %
\]
preserves the measure
$\mu$. Since the Markov chain is irreducible and null recurrent, $T$
is conservative and ergodic (see \cite{harrisrobbins1953}).

Consider the set $A = \{ x \in\bbs^{\bbn}: x_0=i_0 \}$ with the
fixed state $i_0\in\bbs$ chosen above. Let
\[
\varphi_A(x) = \min\bigl\{ n \geq1: 
T^n x
\in A \bigr\},\qquad x \in\bbs^{\bbn} %
\]
be the first entrance time, and assume that
\[
\sum_{k=1}^n P_{i_0}(
\varphi_A \geq k) \in RV_{\beta}, %
\]
the set of regularly varying sequences with exponent $\beta$ of regular
variation, for $\beta\in(0, 1)$. By the Tauberian theorem for power
series (see, e.g., \cite{feller1966}), this is equivalent to assuming
that
%
%
\begin{equation}
\label{ereturnDK} P_{i_0}(\varphi_A \geq k) \in
RV_{\beta-1}.
\end{equation}
There are many natural examples of
Markov chains with this property. Probably, the simplest example is
obtained by taking $\bbs=\{ 0,1,2,\ldots\}$ and letting the transition
probabilities satisfy $p_{i,i-1}=1$ for $i\geq1$, with $ (
p_{0,j},   j=0,1,2,\ldots )$ being an arbitrary probability
distribution satisfying
\[
\sum_{j=k}^\infty p_{0,j}\in
RV_{\beta-1},\qquad k\to\infty. %
\]

Let $f\in L^\infty(\mu)$ be a nonnegative function on $\bbs^{\bbn}$
supported by $A$. Define for $0<\alpha<2$
%
%
\begin{equation}
\label{ebn} b_n = \biggl( \int_E \max
_{1 \leq k \leq n} \bigl( f\circ T^k(x) \bigr)^{\alpha}
\mu(dx) \biggr)^{1/\alpha},\qquad n=1,2,\dots.
\end{equation}
The sequence $(b_n)$ plays an important part in
\cite{owadasamorodnitsky2014}, and it will play an important role in
this paper as well. If we define \textit{the wandering rate sequence} by
\[
w_n = \mu \bigl( \bigl\{ x \in\bbs^{\bbn}:
x_j=i_0\mbox { for some }j=0,1,\ldots, n \bigr\}
\bigr),\qquad n=1,2,\dots, %
\]
then, clearly, $w_n\sim\mu(\varphi_A\leq n)$ as $n\to\infty$. We know
by Theorem 4.1 in \cite{owadasamorodnitsky2014} that
%
%
\begin{equation}
\label{eRVexpbn} \lim_{n\to\infty}\frac{b_n^\alpha}{ w_n}=\llVert f\rrVert
_\infty.
\end{equation}
Furthermore, it
follows from Lemma 3.3 in
\cite{resnicksamorodnitskyxue2000} that
\[
w_n \sim\sum_{k=1}^n
P_{i_0}(\varphi_A \geq k) \in RV_{\beta}.
\]

The above setup allows us to define a stationary symmetric
$\alpha$-stable (S$\alpha$S) sequence by
%
%
\begin{equation}
\label{eunderlyingproc} X_n = \int_E f \circ
T^n(x) \,dM(x),\qquad n=1,2,\dots,
\end{equation}
where $M$ is a S$\alpha$S random measure on
$(E,\mathcal{E})$ with control measure $\mu$. See
\cite{samorodnitskytaqqu1994} for details on $\alpha$-stable random
measures and integrals with respect to these measures. This is a long
range dependent sequence, and the parameter $\beta$ of the Markov
chain determined just how long the memory is; see
\cite{owadasamorodnitsky2015,owadasamorodnitsky2014}. Section~\ref{sec5} of the present paper
discusses an extremal limit theorem for this sequence.

\section{Random sup measure structure}
\label{secsupmeasure}

In this section, we prove a limit theorem, and the limit in this
theorem is a stationary and
self-similar random sup measure whose restrictions to the intervals of
the type $(0,t]$, $t\geq0$, as in~(\ref{erestrM}) is
distributionally equal to the process $Z_{\alpha, \beta}$ in
(\ref{elimprocess}). This result is also a major step toward the
extension of the main result in \cite{owadasamorodnitsky2014} to the
setup in~(\ref{esupmconv}) of weak convergence in the space of sup
measures of normalized partial maxima of the
absolute values of a S$\alpha$S sequence. The extension
itself is formally proved in the next section. We emphasize that the
discussion in this section applies to all $0<\beta<1$.

We introduce first some additional setup. Let $L_{1-\beta}$ be the
standard $(1-\beta)$-stable subordinator,
that is, an increasing L\'evy process such that
\[
Ee^{-\theta L_{1-\beta}(t)} = e^{-t\theta^{1-\beta}}\qquad\mbox {for }\theta\geq0\mbox{ and }t
\geq0. %
\]
Let
%
%
\begin{equation}
\label{erange} R_\beta=\overline{ \bigl\{ L_{1-\beta}(t), t\geq0
\bigr\} }\subset [0,\infty)
\end{equation}
be (the closure of) the range of the subordinator. It has several very
attractive properties as a random closed set, described in the following
proposition. We equip the space $\mathscr J$ of closed subsets of
$[0,\infty)$ with
the usual Fell topology (see \cite{molchanov2005}), and the Borel
$\sigma$-field generated by that topology. We will use some basic
facts about measurability of $\mathscr J$-valued maps and equality of
measures on $\mathscr J$; these are stated in the proof of the
proposition below. It is always sufficient to consider ``hitting''
open sets, and among the latter it is sufficient to consider finite
unions of open intervals.

%
\begin{proposition} \label{prrange}
Let $\beta\in(0,1)$ and $R_{\beta}$ be the range~(\ref{erange}) of
the standard $(1-\beta)$-stable subordinator $L_{1-\beta}$ defined on
some probability space
$ ( \Omega, {\mathcal F}, P )$. Then:
\begin{longlist}[(a)]
\item[(a)]  $R_\beta$ is a random closed subset of $[0,\infty)$.

\item[(b)]  For any $a>0$, $aR_\beta\stackrel{d}{=}R_\beta$ as random
closed sets.

\item[(c)]  Let $\mu_\beta$ be a measure on $(0,\infty)$ given by
$\mu_\beta(dx) = \beta x^{\beta-1}  \,dx,   x>0$, and let
$\kappa_\beta=
(\mu_\beta\times P)\circ H^{-1}$, where $H:  (0,\infty) \times
\Omega\to{\mathscr J}$ is defined by
$H(x,\omega)=R_\beta(\omega)+x$. Then for any $r>0$ the measure
$\kappa_\beta$ is
invariant under the shift map $G_r:  {\mathscr J}\to{\mathscr J}$
given by
\[
G_r(F) =F\cap[r,\infty)-r. %
\]
\end{longlist}
\end{proposition}

\begin{pf}
For part (a), we need to check that for any open $G\subseteq
[0,\infty)$, the set
\[
\bigl\{ \omega\in\Omega: R_\beta(\omega)\cap G\neq\varnothing \bigr\}
\]
is in $\mathcal F$. By the right continuity of sample paths of the
subordinator, the same set can be written in the form
\[
\bigl\{ \omega\in\Omega: L_{1-\beta}(r)\in G\mbox{ for some rational
}r \bigr\}. %
\]
Now the measurability is obvious.

Part (b) is a consequence of the self-similarity of the
subordinator. Indeed, it is enough to check that for any open
$G\subseteq [0,\infty)$
\[
P ( R_\beta\cap G\neq\varnothing ) = P ( aR_\beta\cap G\neq
\varnothing ). %
\]
However, by the self-similarity,
\begin{eqnarray*}
P ( R_\beta\cap G\neq\varnothing ) &=& P \bigl( L_{1-\beta}(r)\in G
\mbox{ for some rational }r \bigr) %
\\
&=& P \bigl( aL_{1-\beta}\bigl(a^{-(1-\beta)}r\bigr)\in G\mbox{ for
some rational }r \bigr) = P ( aR_\beta\cap G\neq\varnothing ),
\end{eqnarray*}
as required.

For part (c) it is enough to check that for any finite collection of
disjoint intervals, $0<b_1<c_1<b_2<c_2<\cdots<b_n<c_n<\infty$
%
%
\begin{eqnarray}\label{echeckshiftm}
&& \kappa_\beta \Biggl( \Biggl\{ F\in{\mathscr J}: F\cap
\bigcup_{j=1}^n (b_j,c_j)
\neq\varnothing \Biggr\} \Biggr)
\nonumber\\[-8pt]\\[-8pt]\nonumber
&&\quad = \kappa_\beta \Biggl( \Biggl\{ F\in{
\mathscr J}: F\cap \bigcup_{j=1}^n
(b_j+r,c_j+r)\neq\varnothing \Biggr\} \Biggr);
\end{eqnarray}
see Example 1.29 in \cite{molchanov2005}. A simple inductive argument
together with the strong Markov property of the subordinator shows
that it is enough to prove~(\ref{echeckshiftm}) for the case of a
single interval. That is, one has to check that for any $0<b<c<\infty$,
%
%
\begin{equation}
\label{echeckshift} \kappa_\beta \bigl( \bigl\{ F\in{\mathscr J}: F\cap
(b,c)\neq\varnothing \bigr\} \bigr) = \kappa_\beta \bigl( \bigl\{ F\in{
\mathscr J}: F\cap (b+r,c+r)\neq\varnothing \bigr\} \bigr).
\end{equation}
For $h>0$, let
\[
\delta_h = \inf \bigl\{ y: y\in R_\beta\cap[h,\infty)
\bigr\} -h %
\]
be the overshoot of the level $h$ by the subordinator
$L_{1-\beta}$. Then~(\ref{echeckshift}) can be restated in the form
\begin{eqnarray*}
&&\int_0^b \beta x^{\beta-1} P(
\delta_{b-x}<c-b) \,dx + \bigl(c^\beta -b^\beta\bigr)
\\
&&\quad  =
\int_0^{b+r} \beta x^{\beta-1} P(
\delta_{b+r-x}<c-b) \,dx + \bigl((c+r)^\beta-(b+r)^\beta
\bigr). %
\end{eqnarray*}
The overshoot $\delta_h$ is known to have a density with respect to
the Lebesgue measure, given by
%
%
\begin{equation}
\label{eovershootd} f_h(y) = \frac{\sin (\pi(1-\beta) )}{\pi} h^{1-\beta}
(y+h)^{-1}y^{\beta-1},\qquad y>0;
\end{equation}
see, for example, Exercise 5.6 in \cite{kyprianou2006}, and checking
the required
identity is a matter of somewhat tedious but still elementary
calculations.
\end{pf}

In the notation of Section~\ref{secergodic}, we define for
$n=1,2,\ldots$ and $x\in E=\bbs^\bbn$ a sup measure on $[0,\infty)$ by
%
%
\begin{equation}
\label{emeanmn} m_n(B; x) = \max_{k:  k/n\in B} f\circ
T^k(x),\qquad B\subseteq[0,\infty ).
\end{equation}

The main result of this section will be stated in terms of weak
convergence of a sequence of finite-dimensional random vectors. Its
significance will go well beyond that weak convergence, as we will
describe in the sequel. Let $0\leq t_1<t_1^\prime
\leq\cdots\leq t_m<t_m^\prime<\infty$ be fixed points, $m\geq1$.
For $n=1,2,\ldots$ let $Y^{(n)}=(Y^{(n)}_1,\ldots, Y^{(n)}_m)$ be an
$m$-dimensional Fr\'echet random vector satisfying
%
%
\begin{equation}
\label{etimen} P \bigl( Y^{(n)}_1\leq\lambda_1,
\ldots, Y^{(n)}_m\leq\lambda _m \bigr) = \exp
\Biggl\{ -\int_E \bigvee_{i=1}^m
\lambda_i^{-\alpha} m_n \bigl(
\bigl(t_i,t_i^\prime\bigr);x
\bigr)^\alpha \mu(dx) \Biggr\},
\end{equation}
for $\lambda_j>0,   j=1,\ldots, m$; see, for example, \cite
{stoevtaqqu2005} for
details on Fr\'echet random vectors and processes.
%
\begin{theorem} \label{tfrechetconv}
Let $0<\beta<1$. The sequence of random vectors $(b_n^{-1}Y^{(n)})$
converges weakly in
$\reals^m$ to a Fr\'echet random vector $Y^\ast=(Y^\ast_1,\ldots,
Y^\ast_m)$ such that
%
%
\begin{eqnarray}\label{etimeinfty}
&& P \bigl( Y^\ast_1\leq\lambda_1,
\ldots, Y^\ast_m\leq\lambda _m \bigr)
\nonumber\\[-8pt]\\[-8pt]\nonumber
&&\quad = \exp
\Biggl\{ - E^\prime \Biggl(\int_0^\infty
\bigvee_{i=1}^m \lambda_i^{-\alpha}
\one \bigl( ( R_\beta+x )\cap\bigl(t_i,t_i^\prime
\bigr)\neq \varnothing \bigr) \beta x^{\beta-1} \,dx \Biggr) \Biggr\}
\end{eqnarray}
for $\lambda_j>0,   j=1,\ldots, m$, where $R_\beta$ is the range
(\ref{erange}) of a $(1-\beta)$-stable subordinator defined on some
probability
space $ ( \Omega^\prime, {\mathcal F}^\prime, P^\prime )$.
\end{theorem}

We postpone proving the theorem and discuss first its
significance. Define
%
%
\begin{equation}
\label{esupmasint}
\hspace*{-10pt}W_{\alpha,\beta}(A) = \hspace*{4pt}\rule{0pt}{13pt}^e\hspace*{-4pt} \int_{(0,\infty)\times
\Omega^\prime}
\one \bigl( \bigl( R_\beta\bigl(\omega^\prime\bigr)+x \bigr)\cap
A\neq\varnothing \bigr) M\bigl(dx,d\omega^\prime\bigr),\qquad A\subseteq [0, \infty), \mbox{ Borel.}
\end{equation}
The integral in~(\ref{esupmasint}) is the extremal integral with
respect to a Fr\'echet random sup measure $M$ on $(0,\infty)\times
\Omega^\prime$, where $ ( \Omega^\prime, {\mathcal F}^\prime,
P^\prime )$ is some probability space. We refer the reader to
\cite{stoevtaqqu2005} for details. The control measure of $M$ is
$m= \mu_\beta\times P^\prime$, where $\mu_\beta$ is defined in
part (c) of
Proposition~\ref{prrange}. It is evident that
$W_{\alpha,\beta}(A)<\infty$ a.s. for any bounded Borel set $A$.
We claim that a version of $W_{\alpha,\beta}$ is a random sup measure on
$[0,\infty)$. The necessity of taking an appropriate version stems
from the usual phenomenon, that
the extremal integral is defined separately for each set $A$, with a
corresponding $A$-dependent exceptional set.

Let $N_{\alpha,\beta}$ be a Poisson random measure on $(0,\infty)^2$
with the mean measure
\[
\alpha x^{-(\alpha+1)} \,dx\, \beta y^{\beta-1} \,dy,\qquad x, y>0. %
\]
Let $ ((U_i,V_i) )$ be a measurable enumeration of the points of
$N_{\alpha,\beta}$. Let, further, $ ( R_\beta^{(i)} )$ be
i.i.d. copies of the range of the $(1-\beta)$-stable
subordinator, independent of the Poisson random measure
$N_{\alpha,\beta}$. Then a version of $W_{\alpha,\beta}$ is given by
%
%
\begin{equation}
\label{ewtilde} \hat W_{\alpha,\beta}(A) = \bigvee_{i=1}^\infty
U_i \one \bigl( \bigl( R_\beta^{(i)}+V_i
\bigr)\cap A\neq\varnothing \bigr),\qquad A\subseteq [0,\infty),\mbox{ Borel;}
\end{equation}
see \cite{stoevtaqqu2005}. It is interesting to note that, since
the origin belongs, with probability 1, to the range of the
subordinator, evaluating~(\ref{ewtilde}) on sets of the form
$A=[0,t]$, $0\leq t\leq1$, reduces this representation to the more
standard representation of the process $Z_{\alpha, \beta}$ in
(\ref{elimprocess}). See (3.8) in \cite{owadasamorodnitsky2014}.

It is clear that $\hat W_{\alpha,\beta}$
is a random sup measure on $[0,\infty)$. In fact,
%
%
\begin{equation}
\label{eWcheck} d\, \check{\phantom{_i}}\, \hat W_{\alpha,\beta}(t) = \cases{ U_i, &
\quad if $t\in R_\beta^{(i)}+V_i$, some $i=1,2,
\ldots,$
\vspace*{3pt}\cr
0, &\quad otherwise.}
\end{equation}
Even though it is $\hat W_{\alpha,\beta}$ that takes values in
the space of sup measures, we will slightly abuse the terminology and
refer to $W_{\alpha,\beta}$ itself a random sup measure.

%
\begin{proposition} \label{prstatssrsm}
For any $\beta\in(0,1)$, the random sup measure $W_{\alpha,\beta}$
is stationary and
self-similar with exponent $H=\beta/\alpha$ in the sense of
(\ref{esasupm}).
\end{proposition}

\begin{pf}
Both\vspace*{1pt} statements can be read off~(\ref{eWcheck}). Indeed, the pairs
$ (U_i,(R_\beta^{(i)}+V_i) )$ form a Poisson random measure on
$(0,\infty)\times{\mathscr J}$ and, by part (c) of Proposition
\ref{prrange}, the mean measure of this Poisson random measure is
unaffected by the transformations $G_r$ applied to the random set
dimension. This implies the law of the random upper semicontinuous
function $d \, \check{\phantom{_i}}  \hat W_{\alpha,\beta}$ is shift invariant,
hence stationarity of $W_{\alpha,\beta}$.

For the self-similarity, note that replacing $t$ by $t/a$, $a>0$ in
(\ref{eWcheck}) is equivalent to replacing $R_\beta^{(i)}$ by
$aR_\beta^{(i)}$ and $V_i$ by $aV_i$. By part (b) of Proposition
\ref{prrange}, the former action does not change the law of a
random closed set, while it is elementary to check that the law of the
Poisson random measure on $(0,\infty)^2$ with points
$ ((U_i,aV_i) )$ is the same as the law of the Poisson random
measure on the same space with the points
$ ((a^{\beta/\alpha}U_i,V_i) )$. Hence, the self-similarity
of $W_{\alpha,\beta}$ with $H=\beta/\alpha$.
\end{pf}

Returning now to the result in Theorem~\ref{tfrechetconv}, note that it
can be restated in the form
\[
b_n^{-1}\bigl(Y^{(n)}_1,\ldots,
Y^{(n)}_m\bigr) \Rightarrow \bigl( W_{\alpha,\beta}\bigl(
\bigl(t_1,t_1^\prime\bigr)\bigr), \ldots,
W_{\alpha,\beta}\bigl(\bigl(t_m,t_m^\prime\bigr)
\bigr) \bigr)\qquad\mbox{as }n\to \infty. %
\]
In particular, if we choose $t_i =t_{i-1}^\prime,   i=1,\ldots, m$,
with $t_1=0$ and an arbitrary $t_{m+1}$, and define
\[
Z^{(n)}_i=\max_{j=1,\ldots,i}Y^{(n)}_j,\qquad
i=1,\ldots, m, %
\]
then
%
%
\begin{eqnarray}
\label{eytomeas} \bigl(b_n^{-1}Z^{(n)}_i,
i=1,\ldots, m\bigr) &\Rightarrow& \Bigl( \max_{j=1,\ldots,i}W_{\alpha,\beta}
\bigl((t_j,t_{j+1})\bigr), i=1,\ldots, m \Bigr)
\nonumber\\[-8pt]\\[-8pt]\nonumber
&  =& \bigl( W_{\alpha,\beta}\bigl((0,t_{i+1})
\bigr), i=1,\ldots, m \bigr).
\nonumber
\end{eqnarray}
However, as a part of the argument in \cite{owadasamorodnitsky2014}
it was established that
\[
\bigl(b_n^{-1}Z^{(n)}_i, i=1,\ldots,
m\bigr) \Rightarrow \bigl( Z_{\alpha,\beta}(t_{i+1}), i=1,\ldots, m
\bigr), %
\]
with $Z_{\alpha,\beta}$ as in~(\ref{elimprocess}); this is (4.7) in \cite{owadasamorodnitsky2014}.
This leads to the
immediate conclusion, stated in the following corollary.
%
\begin{corollary} \label{cmainresult}
For any $\beta\in(0,1)$, the time-changed extremal Fr\'echet process satisfies
\[
\bigl( Z_{\alpha,\beta}(t), t\geq0 \bigr) \stackrel{d} {=} \bigl( W_{\alpha,\beta}
\bigl((0,t]\bigr), t\geq0 \bigr) %
\]
and, hence, is a restriction of the stationary and self-similar random
sup measure $W_{\alpha,\beta}$ (to the intervals $(0,t],   t\geq0$).
\end{corollary}

We continue with a preliminary result, needed for the proof of Theorem
\ref{tfrechetconv}, which may also be of independent
interest.
%
\begin{proposition} \label{prCLTextreme}
Let $0<\gamma<1$, and $(Y_1,Y_2,\ldots)$ be i.i.d. nonnegative random
variables such that $P(Y_1>y)$ is regularly varying with exponent
$-\gamma$. Let $S_0=0$ and $S_n=Y_1+\cdots+Y_n$ for $n=1,2,\ldots$ be
the corresponding partial sums. For $\theta>0$ define a random sup
measure on $[0,\infty)$ by
\[
M^{(Y;\theta)}(G) = \one ( S_n\in\theta G\mbox{ for some
}n=0,1,\ldots ) %
\]
$G\subseteq[0,\infty)$, open. Then
\[
M^{(Y;\theta)}\Rightarrow_{\theta\to\infty} M^{(\gamma)} %
\]
in the space SM equipped with the sup vague topology, where
\[
M^{(\gamma)}(G) = \one ( R_{1-\gamma}\cap G\neq\varnothing ). %
\]
\end{proposition}

\begin{pf}
It is enough to prove that for any finite collection of intervals
$(a_i,b_i), i=1,\ldots, m$ with $0<a_i<b_i<\infty,   i=1,\ldots, m$
we have
%
%
\begin{eqnarray}\label{eCLTextrint}
&& P \bigl(\mbox{for each }i=1,\ldots, m, S_j/\theta\in
(a_i,b_i)\mbox{ for some }j=1,2,\ldots \bigr)
\nonumber\\[-8pt]\\[-8pt]\nonumber
&&\quad \to P \bigl(\mbox{for each }i=1,\ldots, m, R_{1-\gamma}\cap
(a_i,b_i)\neq\varnothing \bigr) %
\end{eqnarray}
as $\theta\to\infty$. If we let $a(\theta) =  (
P(Y_1>\theta) )^{-1}$, a regularly varying function with exponent
$\gamma$, then the probability in the left-hand side
of~(\ref{eCLTextrint}) can be rewritten as
%
%
\begin{equation}
\label{eCLTtint} P \bigl( \mbox{for each }i=1,\ldots, m, S_{\lfloor
ta(\theta)\rfloor}/\theta
\in(a_i,b_i) \mbox{ for some }t\geq0 \bigr).
\end{equation}
By the invariance principle,
%
%
\begin{equation}
\label{einvpr} ( S_{\lfloor ta(\theta)\rfloor}/\theta, t\geq0 ) \Rightarrow_{\theta\to\infty}
\bigl( L_\gamma(t), t\geq0 \bigr)
\end{equation}
weakly in the $J_1$-topology in the space $D[0,\infty)$, where
$L_\gamma$ is the standard $\gamma$-stable subordinator; see,
for example, \cite{jacodshiryaev1987}. If we denote by
$D_+^\uparrow[0,\infty)$ the set of all nonnegative nondecreasing
functions in $D[0,\infty)$ vanishing at $t=0$, then
$D_+^\uparrow[0,\infty)$ is, clearly,
a closed set in the $J_1$-topology, so the weak convergence in
(\ref{einvpr}) also takes places in the $J_1$-topology relativized
to $D_+^\uparrow[0,\infty)$.

For a function $\varphi\in D_+^\uparrow[0,\infty)$, let
\[
R_\varphi= \overline{ \bigl\{ \varphi(t), t\geq0 \bigr\}} %
\]
be the closure of its range. Notice that
\[
R_\varphi= \biggl( \bigcup_{t> 0} \bigl(
\varphi(t-),\varphi(t) \bigr) \biggr)^c, %
\]
which makes it evident that for any $0<a<b<\infty$ the set
\[
\bigl\{ \varphi\in D_+^\uparrow[0,\infty): R_\varphi\cap [a,b]=
\varnothing \bigr\} %
\]
is open in the $J_1$-topology, hence measurable. Therefore, the set
\[
\bigl\{ \varphi\in D_+^\uparrow[0,\infty): R_\varphi\cap (a,b)\neq
\varnothing \bigr\} = \bigcup_{k=1}^\infty \bigl
\{ \varphi\in D_+^\uparrow[0,\infty): R_\varphi\cap[a+1/k,b-1/k]\neq
\varnothing \bigr\} %
\]
is measurable as well and, hence, so is the set
\[
\bigl\{ \varphi\in D_+^\uparrow[0,\infty): \mbox{for each }i=1,\ldots, m,
R_\varphi\cap (a_i,b_i)\neq\varnothing \bigr\}.
\]
Therefore, the desired conclusion~(\ref{eCLTextrint}) will follow
from~(\ref{eCLTtint}) and the invariance principle~(\ref{einvpr})
once we check that the measurable function on $D_+^\uparrow[0,\infty)$
defined by
\[
J(\varphi) = \one \bigl( R_\varphi\cap (a_i,b_i)
\neq\varnothing\mbox{ for each }i=1,\ldots, m \bigr) %
\]
is a.s. continuous with respect to the law of $L_\gamma$ on
$D_+^\uparrow[0,\infty)$. To see this, let
\[
B_1 = \bigl\{ \varphi\in D_+^\uparrow[0,\infty): \mbox{for
each }i=1,\ldots, m\mbox{ there is } t_i\mbox{ such that }
\varphi(t_i)\in (a_i,b_i) \bigr\}
\]
and
\begin{eqnarray*}
B_2 &=& \bigl\{ \varphi\in D_+^\uparrow[0,\infty): \mbox{for some
}i=1,\ldots, m\mbox{ there is }t_i
\mbox{ such that }(a_i,b_i)\subseteq \bigl(
\varphi(t_i-), \varphi (t_i) \bigr) \bigr\}. %
\end{eqnarray*}
Both sets are open in the $J_1$-topology on $D_+^\uparrow[0,\infty)$,
and $J(\varphi)=1$ on $B_1$ and $J(\varphi)=0$ on $B_2$. Now the
a.s. continuity of the function $J$ follows from the fact that
\[
P(L_\gamma\in B_1\cup B_2)=1, %
\]
since a stable subordinator does not hit fixed points.
\end{pf}

%
\begin{remark} \label{rkalsosets}
It follows immediately from Proposition~\ref{prCLTextreme} that we
also have weak convergence in the space of closed subsets of
$[0,\infty)$. Specifically, the random closed set
$\theta^{-1} \{ S_n,  n=0,1,\ldots\}$ converges weakly, as
$\theta\to\infty$, to the random closed set $R_{1-\gamma}$.
\end{remark}

\begin{pf*}{Proof of Theorem~\ref{tfrechetconv}}
We will prove that
%
%
\begin{eqnarray}\label{eallint}
&& \frac{\int_E \min_{i=1,\ldots, m} m_n (
(t_i,t_i^\prime);x )^\alpha  \mu(dx)}{\int_E \max_{1 \leq k
\leq n}  ( f\circ T^n(x) )^{\alpha} \mu(dx)}
\nonumber\\[-8pt]\\[-8pt]\nonumber
&&\quad\to\int_0^\infty\beta x^{\beta-1}
P^\prime \bigl( ( R_\beta +x )\cap \bigl(t_i,t_i^\prime
\bigr)\neq\varnothing\mbox{ for each }i=1,\ldots, m \bigr) \,dx %
\end{eqnarray}
as $n\to\infty$. The reason this will suffice for the proof of the
theorem is that, by the inclusion--exclusion formula, the expression in
the exponent in the right-hand side of~(\ref{etimeinfty}) can be
written as a finite linear combination of terms of the form of the
right-hand side of~(\ref{eallint}) (with different collections of
intervals in each term). More specifically, we can write, for a fixed
$x>0$,
\begin{eqnarray*}
&& E^\prime \Biggl(\bigvee_{i=1}^m
\lambda_i^{-\alpha} \one \bigl( ( R_\beta+x )\cap
\bigl(t_i,t_i^\prime\bigr)\neq\varnothing \bigr)
\Biggr)
\\
&&\quad
= \int_0^\infty P^\prime \bigl(
(R_\beta+x )\cap \bigl(t_i,t_i^\prime
\bigr)\neq\varnothing\mbox{ for some $i$ such that } \lambda_i^{-\alpha}>u
\bigr) \,du %
\end{eqnarray*}
and\vspace*{1pt} apply the inclusion--exclusion formula to the probability of the
union inside the integral. A~similar relation exists between the
left-hand side of~(\ref{eallint}) and the distribution of
$(b_n^{-1}Y^{(n)})$.

An additional simplification that we may and will introduce is that of
assuming that $f$ is constant on $A$. Indeed, it follows immediately
from the ergodicity that both the numerator and the denominator in the
left-hand side of~(\ref{eallint}) do not change asymptotically if
we replace $f$ by $\llVert   f\rrVert  _\infty\one_A$; see (4.2) in
\cite{owadasamorodnitsky2014}. With this simplification,
(\ref{eallint}) reduces to the following statement: as
$n\to\infty$,
%
%
\begin{eqnarray}\label{eallintsimple}
&& \frac{1}{w_n}\mu \Biggl(\bigcap_{i=1}^m
\bigl\{ x_k=i_0\mbox{ for some $k$ with }
t_i<k/n<t_i^\prime \bigr\} \Biggr)
\nonumber\\[-8pt]\\[-8pt]\nonumber
&&\quad \to \int_0^\infty\beta x^{\beta-1}
P^\prime \bigl( ( R_\beta +x )\cap \bigl(t_i,t_i^\prime
\bigr)\neq\varnothing\mbox{ for each }i=1,\ldots, m \bigr) \,dx. %
\end{eqnarray}
Note that we have used~(\ref{eRVexpbn}) in translating
(\ref{eallint}) into the form~(\ref{eallintsimple}).

We introduce the notation $A_0=A$, $A_k = A^c \cap\{\varphi_A = k \}$
for $k \geq
1$. Let $(Y_1,Y_2,\ldots)$ be a sequence of i.i.d. $\bbn$-valued
random variables defined on some probability
space $ ( \Omega^\prime, {\mathcal F}^\prime, P^\prime )$
such that
$P^\prime(Y_1=k)=P_{i_0}(\varphi_A=k)$,  $k=1,2,\ldots.$ By our
assumption, the probability tail $P(Y_1>y)$ is regularly varying with
exponent $-(1-\beta)$. With $S_0=0$ and
$S_j=Y_1+\cdots+Y_j$ for $j=1,2,\ldots$ we have
\begin{eqnarray*}
&& \mu \Biggl(\bigcap_{i=1}^m \bigl\{
x_k=i_0\mbox{ for some }k\mbox{ with }
t_i<k/n<t_i^\prime \bigr\} \Biggr)
\\
&&\quad = \sum_{l:  l/n\leq t_1} \mu(A_l)P^\prime
\bigl( \mbox{for each }i=1,\ldots, m, S_j\in\bigl(nt_i-l,
nt_i^\prime-l\bigr)\\
&&\qquad{} \mbox{for some }j=0,1,\ldots \bigr)
\\
&&\qquad{} + \sum_{l:  t_1<l/n< t_1^\prime} \mu(A_l)P^\prime
\bigl( \mbox{for each }i=2,\ldots, m, S_j\in\bigl(nt_i-l,
nt_i^\prime-l\bigr) \mbox{ for some }j=0,1,\ldots
\bigr)
\\
&&\quad := D_n^{(1)}+D_n^{(2)}. %
\end{eqnarray*}
It is enough to prove that
%
%
\begin{equation}
\label{erange1} \lim_{n\to\infty} \frac{1}{w_n}
D_n^{(1)}= \int_0^{t_1}
\beta x^{\beta-1} P^\prime \bigl( ( R_\beta+x )\cap
\bigl(t_i,t_i^\prime\bigr)\neq\varnothing
\mbox{ for each }i=1,\ldots, m \bigr) \,dx
\end{equation}
and
%
%
\begin{equation}
\label{erange2} \lim_{n\to\infty} \frac{1}{w_n}D_n^{(2)}
= \int_{t_1}^{t_1^\prime
} \beta x^{\beta-1}
P^\prime \bigl( ( R_\beta+x )\cap \bigl(t_i,t_i^\prime
\bigr)\neq\varnothing\mbox{ for each }i=1,\ldots, m \bigr) \,dx.
\end{equation}
We will prove~(\ref{erange1}), and~(\ref{erange2}) can be proved
in the same way. Let $K$ be a large positive integer, and
$\varepsilon>0$ a small number. For each integer $1\leq d\leq
(1-\epsilon) K$, and each $l:  t_1(d-1)/K \leq l/n<t_1d/K$, we have
\begin{eqnarray*}
\hspace*{-4pt}&& P^\prime \bigl( \mbox{for each }i=1,\ldots, m, S_j\in
\bigl(nt_i-l, nt_i^\prime-l\bigr)
\mbox{ for some }j=0,1,\ldots \bigr) %
\\
\hspace*{-4pt}&&\quad \leq P^\prime \bigl( \mbox{for each }i=1,\ldots, m, S_j\in
\bigl(nt_i-nt_1d/K, nt_i^\prime-nt_1(d-1)/K
\bigr)\\
\hspace*{-4pt}&&\qquad{} \mbox{for some } j=0,1,\ldots \bigr) %
\\
\hspace*{-4pt}&&\quad \to P^\prime \bigl(\mbox{for each }i=1,\ldots, m, R_{\beta}\cap
\bigl( t_i-t_1d/K, t_i^\prime-t_1(d-1)/K
\bigr)\neq\varnothing \bigr) %
\end{eqnarray*}
as $n\to\infty$, by Proposition~\ref{prCLTextreme}. Therefore,
\begin{eqnarray*}
&& \limsup_{n\to\infty} \frac{1}{w_n} D_n^{(1)}
\\
&&\quad
\leq\sum_{d=1}^{\lfloor
(1-\epsilon) K \rfloor} \biggl[\limsup
_{n\to\infty} \frac{\sum_{l:
t_1(d-1)/K \leq l/n<t_1d/K} \mu(A_l)}{w_n}
\\
&&\qquad {}\times P^\prime \bigl( \mbox{for each }i=1,\ldots, m, R_{\beta}\cap
\bigl( t_i-t_1d/K, t_i^\prime-t_1(d-1)/K
\bigr)\neq\varnothing \bigr) \biggr]
\\
&&\qquad{} + \limsup_{n\to\infty} \frac{\sum_{l:
t_1 \lfloor
(1-\epsilon) K \rfloor/K\leq l/n\leq t_1 } \mu(A_l)}{w_n}. %
\end{eqnarray*}
Since for any $a>0$,
\[
\sum_{l=1}^{na}\mu(A_l) \sim
w_{\lfloor na\rfloor} \qquad\mbox{as }n\to\infty, %
\]
and the wandering sequence $(w_n)$ is regularly varying with exponent
$\beta$, we conclude that
\begin{eqnarray*}
\limsup_{n\to\infty} \frac{\sum_{l:
t_1(d-1)/K \leq l/n<t_1d/K} \mu(A_l)}{w_n} &=& \limsup_{n\to\infty}
\frac{ w_{\lfloor nt_1d/K\rfloor} -
w_{\lfloor
nt_1(d-1)/K\rfloor} }{w_n} %
\\
&=& \frac{t_1^\beta}{K^\beta} \bigl( d^\beta- (d-1)^\beta \bigr)
\end{eqnarray*}
for $1\leq d\leq(1-\epsilon) K$ and, similarly,
\[
\limsup_{n\to\infty} \frac{\sum_{l:
t_1 \lfloor
(1-\epsilon) K \rfloor/K\leq l/n\leq t_1 } \mu(A_l)}{w_n} = t_1^\beta
\biggl[ 1- \biggl( \frac{\lfloor(1-\varepsilon)K\rfloor
}{K} \biggr)^\beta \biggr].
\]
Therefore,
\begin{eqnarray*}
&& \limsup_{n\to\infty} \frac{1}{w_n} D_n^{(1)}
\\
&&\quad \leq\int_0^{(1-\varepsilon)t_1} \beta x^{\beta-1}
P^\prime \bigl( R_\beta\cap \bigl( t_i-a_K(x),
t_i^\prime-b_K(x) \bigr) \neq\varnothing
\mbox{ for each }i=1,\ldots, m \bigr) \,dx %
\\
&&\qquad{} + t_1^\beta \biggl[ 1- \biggl( \frac{\lfloor(1-\varepsilon)K\rfloor
}{K}
\biggr)^\beta \biggr], %
\end{eqnarray*}
where $a_K(x)= t_1d/K$ and $b_K(x) = t_1(d-1)/K$ if $t_1(d-1)/K\leq
x<t_1d/K$ for $1\leq d\leq(1-\epsilon) K$. Since
\[
\one \bigl( R_\beta\cap(a_k,b_k)\neq\varnothing
\bigr)\to \one \bigl( R_\beta\cap(a,b)\neq\varnothing \bigr) %
\]
a.s. if $a_k\to a$ and $b_k\to b$, we can let $K\to\infty$ and then
$\varepsilon\to0$ to conclude that
%
%
\begin{equation}
\label{eupperbmeas} \limsup_{n\to\infty} \frac{1}{w_n}
D_n^{(1)} \leq \int_0^{t_1}
\beta x^{\beta-1} P^\prime \bigl( R_\beta\cap \bigl(
t_i-x, t_i^\prime-x \bigr) \neq\varnothing \mbox{ for each }i=1,\ldots, m \bigr) \,dx.
\end{equation}

We can obtain a lower bound matching~(\ref{eupperbmeas}) in a
similar way. Indeed, for each integer $1\leq d\leq
(1-\epsilon) K$, and each $l:  t_1(d-1)/K \leq l/n<t_1d/K$ as above,
we have
\begin{eqnarray*}
\hspace*{-4pt}&& P^\prime \bigl( \mbox{for each }i=1,\ldots, m, S_j\in
\bigl(nt_i-l, nt_i^\prime-l\bigr) \mbox{ for
some }j=0,1,\ldots \bigr) %
\\
\hspace*{-4pt}&&\quad \geq P^\prime \bigl( \mbox{for each }i=1,\ldots, m, S_j\in
\bigl(nt_i-nt_1(d-1)/K, nt_i^\prime-nt_1d/K
\bigr)\\
\hspace*{-4pt}&&\qquad{} \mbox{for some }j=0,1,\ldots \bigr) %
\\
\hspace*{-4pt}&&\quad \to P^\prime \bigl( \mbox{for each }i=1,\ldots, m, R_{\beta
}\cap
\bigl( t_i-t_1(d-1)/K, t_i^\prime-t_1d/K
\bigr)\neq\varnothing \bigr) %
\end{eqnarray*}
as $n\to\infty$, by Proposition~\ref{prCLTextreme}, and we proceed
as before. This gives a lower bound complementing
(\ref{eupperbmeas}), so we have proved that
\[
\lim_{n\to\infty} \frac{1}{w_n} D_n^{(1)} =
\int_0^{t_1} \beta x^{\beta-1}
P^\prime \bigl( R_\beta\cap \bigl( t_i-x,
t_i^\prime-x \bigr) \neq\varnothing \mbox{ for each }i=1,
\ldots, m \bigr) \,dx. %
\]
This is, of course,~(\ref{erange1}).
\end{pf*}

\section{Convergence in the space SM}
\label{secconvergence}\label{sec5}

Let $\BX=(X_1,X_2,\ldots)$ be the stationary S$\alpha$S process
defined by
(\ref{eunderlyingproc}). The following theorem is a partial
extension of Theorem 4.1 in \cite{owadasamorodnitsky2014} to weak
convergence in the space of sup measures. In its statement, we use the
usual tail constant of an $\alpha$-stable random variable given by
\[
C_{\alpha} = \biggl( \int_0^{\infty}
x^{-\alpha} \sin x \,dx \biggr)^{-1} = \cases{ (1-\alpha) / \bigl(
\Gamma(2-\alpha) \cos(\pi \alpha/ 2) \bigr), &\quad if $\alpha\neq1$,
\cr
2 / \pi,
&\quad if $\alpha=1$;} %
\]
see \cite{samorodnitskytaqqu1994}.

%
\begin{theorem} \label{textendOS}
For $n=1,2,\ldots$ define a random sup measure $M_n(\llvert  \BX\rrvert  )$ on
$[0,\infty)$ by~(\ref{epartialmn}), with
$\llvert  \BX\rrvert  =(\llvert  X_1\rrvert,\llvert  X_2\rrvert,\ldots)$. Let $(b_n)$ be given by~(\ref{ebn}). If
$1/2<\beta<1$, then
%
%
\begin{equation}
\label{eextendOS} \frac{1}{b_n}M_n\bigl(\llvert \BX\rrvert \bigr)
\Rightarrow C_\alpha^{1/\alpha}W_{\alpha,\beta}\qquad\mbox{as }n\to
\infty
\end{equation}
in the sup vague topology in the space SM.
\end{theorem}

\begin{pf}
The weak convergence in the space SM will be established if we show
that for any $0\leq t_1<t_1^\prime
\leq\cdots\leq t_m<t_m^\prime<\infty$,
\begin{eqnarray*}
&&\bigl( b_n^{-1} M_n\bigl(\llvert \BX\rrvert
\bigr) \bigl( \bigl(t_1,t_1^\prime\bigr) \bigr),
\ldots, b_n^{-1} M_n\bigl(\llvert \BX\rrvert
\bigr) \bigl( \bigl(t_m,t_m^\prime\bigr) \bigr)
\bigr) \Rightarrow C_\alpha^{1/\alpha}\bigl( W_{\alpha,\beta} \bigl( \bigl(t_1,t_1^\prime
\bigr) \bigr), \ldots,\\
&&\quad W_{\alpha,\beta} \bigl( \bigl(t_m,t_m^\prime
\bigr) \bigr) \bigr) %
\end{eqnarray*}
as $n\to\infty$ (see Section 12.7 in \cite{vervaat1997}). %
For simplicity
of notation, we will assume that $t_m^\prime\leq1$. Our goal is, then,
to show that
%
%
\begin{equation}
\label{econvSMint} \biggl( \frac{1}{b_n} \max_{n t_1<k<nt_1^\prime}\llvert
X_k\rrvert, \ldots, \frac{1}{b_n} \max_{n t_m<k<nt_m^\prime}
\llvert X_k\rrvert \biggr) \Rightarrow C_\alpha^{1/\alpha}\bigl( W_{\alpha,\beta} \bigl(
\bigl(t_1,t_1^\prime\bigr) \bigr), \ldots,
W_{\alpha,\beta} \bigl( \bigl(t_m,t_m^\prime
\bigr) \bigr) \bigr)
\end{equation}
as $n\to\infty$.

We proceed in the manner similar to that adopted in
\cite{owadasamorodnitsky2014}, and use a series representation of
the S$\alpha$S sequence $(X_1,X_2,\ldots)$. Specifically, we have
%
%
\begin{equation}
\label{eseriesmax} (X_k, k=1,\dots,n) \stackrel{d} {=} \Biggl(
b_n C_{\alpha}^{1/\alpha} \sum
_{j=1}^{\infty} \epsilon_j
\Gamma_j^{-1/\alpha} \frac{f\circ T^k(U_j^{(n)})}{\max_{1 \leq i
\leq n}
f\circ T^i(U_j^{(n)})}, k=1,\dots,n \Biggr).
\end{equation}
In the right-hand side, $(\epsilon_j)$ are i.i.d. Rademacher random
variables
(symmetric random variables with values $\pm
1$), $(\Gamma_j)$ are the arrival times of a unit rate Poisson
process on $(0,\infty)$, and $(U_j^{(n)})$ are i.i.d. $E$-valued random
variables with the common law $\eta_n$ defined by
%
%
\begin{equation}
\label{eetan} \frac{d\eta_n}{d\mu}(x) = \frac{1}{b_n^{\alpha}} \max
_{1 \leq k
\leq
n} f\circ T^k(x)^{\alpha},\qquad x \in E.
\end{equation}
The three sequences $(\epsilon_j)$, $(\Gamma_j)$ and $(U_j^{(n)})$ are
independent. We refer the reader to Section~3.10 of
\cite{samorodnitskytaqqu1994} for details on series representations of
$\alpha$-stable processes. We will prove that for any $\lambda_i>0,
  i=1,\ldots, m$ and $0<\delta<1$,
%
%
\begin{eqnarray}\label{euppbddmax}
&&P \Bigl( b_n^{-1} \max_{nt_i < k < nt_i^\prime} \llvert
X_k\rrvert > \lambda_i,  i=1,\dots,m \Bigr)
\nonumber\\[-8pt]\\[-8pt]\nonumber
&&\quad \leq P \Biggl( C_{\alpha}^{1/\alpha} \bigvee
_{j=1}^{\infty} \Gamma_j^{-1/\alpha}
\frac{\max_{ nt_i < k < nt_i^\prime}f\circ
T^k(U_j^{(n)})}{\max_{1 \leq k \leq
n}f\circ T^k(U_j^{(n)})} > \lambda_i(1-\delta), i=1,\dots, m \Biggr) + o(1)
\end{eqnarray}
and that
%
%
\begin{eqnarray}\label{elowerbddmax}
&&P \Bigl( b_n^{-1} \max_{nt_i < k < nt_i^\prime} \llvert
X_k\rrvert > \lambda_i,  i=1,\dots,m \Bigr)
\nonumber\\[-8pt]\\[-8pt]\nonumber
&&\quad \geq P \Biggl( C_{\alpha}^{1/\alpha} \bigvee
_{j=1}^{\infty} \Gamma_j^{-1/\alpha}
\frac{\max_{ nt_i < k < nt_i^\prime}f\circ
T^k(U_j^{(n)})}{\max_{1 \leq k \leq
n}f\circ T^k(U_j^{(n)})} > \lambda_i(1+\delta), i=1,\dots, m \Biggr) + o(1)
\end{eqnarray}
as $n\to\infty$. Before doing so, we will make a few simple
observations. Let
\[
V_i^{(n)} = \bigvee_{j=1}^{\infty}
\Gamma_j^{-1/\alpha} \frac{\max_{ nt_i < k < nt_i^\prime}f\circ
T^k(U_j^{(n)})}{\max_{1 \leq k \leq
n}f\circ T^k(U_j^{(n)})},\qquad i=1,\ldots, m.
\]
Since the points in $\reals^m$ given by
\[
\biggl( \Gamma_j^{-1/\alpha} \frac{\max_{ nt_i < k < nt_i^\prime
}f\circ
T^k(U_j^{(n)})}{\max_{1 \leq k \leq
n}f\circ T^k(U_j^{(n)})}, i=1,\ldots, m
\biggr),\qquad j=1,2,\ldots %
\]
form a Poisson random measure on $\reals^m$, say, $N_P$, for
$\lambda_i>0$, $ i=1,\ldots,m$ we can write
\begin{eqnarray*}
P \bigl( V^{(n)}_1\leq\lambda_1,\ldots,
V^{(n)}_m\leq\lambda _m \bigr) &=& P \bigl(
N_P \bigl( D(\lambda_1,\ldots, \lambda_m)=0
\bigr) \bigr)
\\
&=& \exp \bigl\{ - E \bigl( N_P \bigl( D(\lambda_1,
\ldots, \lambda_m) \bigr) \bigr) \bigr\},
\end{eqnarray*}
where
\[
D(\lambda_1,\ldots, \lambda_m) = \bigl\{
(z_1,\ldots, z_m): z_i>
\lambda_i \mbox{ for some }i=1,\ldots, m \bigr\}. %
\]
Evaluating the expectation, we conclude that, in the notation of
(\ref{emeanmn}),
\[
P \bigl( V^{(n)}_1\leq\lambda_1,\ldots,
V^{(n)}_m\leq\lambda _m \bigr) = \exp \Biggl\{
-b_n^{-\alpha}\int_E \bigvee
_{i=1}^m \lambda _i^{-\alpha}
m_n \bigl( \bigl(t_i,t_i^\prime
\bigr);x \bigr)^\alpha \mu(dx) \Biggr\}. %
\]
By~(\ref{etimen}), this shows that, in the notation of Theorem
\ref{tfrechetconv},
\[
\bigl( V^{(n)}_1, \ldots, V^{(n)}_m
\bigr) \stackrel{d} {=} \bigl( b_n^{-1}Y^{(n)}_1,
\ldots, b_n^{-1}Y^{(n)}_m \bigr).
\]
Now Theorem~\ref{tfrechetconv} along with the discussion following
the statement of that theorem, and the continuity of the Fr\'echet
distribution show that~(\ref{econvSMint}) and, hence, the claim of
the present theorem, will follow once we prove~(\ref{euppbddmax})
and~(\ref{elowerbddmax}). The two statements can be proved in a
very similar way, so we only prove~(\ref{euppbddmax}).

Once again, we proceed as in \cite{owadasamorodnitsky2014}.
Choose constants $K \in\bbn$ and $0 < \epsilon< 1$ such that both
\[
K+1 > \frac{4}{\alpha} \quad\mbox{and}\quad \delta- \epsilon K > 0.
\]
Then
\begin{eqnarray*}
&& P \Bigl( b_n^{-1} \max_{nt_i < k < nt_i^\prime} \llvert
X_k\rrvert > \lambda _i, i=1,\dots,m \Bigr)
\\
&&\quad \leq P \Biggl( C_{\alpha}^{1/\alpha} \bigvee
_{j=1}^{\infty} \Gamma_j^{-1/\alpha}
\frac{\max_{ nt_i < k < nt_i^\prime}f\circ
T^k(U_j^{(n)})}{\max_{1 \leq k \leq
n}f\circ T^k(U_j^{(n)})} > \lambda_i(1-\delta), i=1,\dots, m \Biggr)
\\
&&\qquad{}+ \varphi_n \Bigl( C_{\alpha}^{-1/\alpha} \epsilon \min
_{1 \leq i \leq m} \lambda_i \Bigr) + \sum
_{i=1}^m \psi _n\bigl(
\lambda_i, t_i,t_i^\prime\bigr),
\end{eqnarray*}
where
\[
\varphi_n(\eta) = P \Biggl( \bigcup_{k=1}^n
\biggl\{ \Gamma_j^{-1/\alpha} \frac{f\circ T^k(U_j^{(n)})}{\max_{1 \leq i
\leq
n}f\circ T^i(U_j^{(n)})} > \eta \mbox{ for at least 2 different } j=1,2,\dots \biggr\} \Biggr), %
\]
$\eta>0$, and for $t<t^\prime$,
\begin{eqnarray*}
\psi_n\bigl(\lambda, t, t^\prime\bigr)
&=& P \Biggl(
C_{\alpha}^{1/\alpha} \max_{nt
< k < nt^\prime} \Biggl\llvert \sum
_{j=1}^{\infty} \epsilon_j
\Gamma_j^{-1/\alpha} \frac{f\circ T^k(U_j^{(n)})}{\max_{1 \leq i
\leq n}f\circ T^i(U_j^{(n)})} \Biggr\rrvert > \lambda,
\\
&&{} C_{\alpha}^{1/\alpha} \bigvee_{j=1}^{\infty}
\Gamma_j^{-1/\alpha} \frac{\max_{nt< k < nt^\prime}f\circ T^k(U_j^{(n)})}{\max_{1 \leq k
\leq n}f\circ T^k(U_j^{(n)})} \leq\lambda(1-\delta),
 \mbox{ and for each } l=1,\dots,n,
\\
&&{} C_{\alpha}^{1/\alpha} \Gamma_j^{-1/\alpha}
\frac{f\circ
T^l(U_j^{(n)})}{\max_{1 \leq i \leq n}f\circ T^i(U_j^{(n)})} > \epsilon\lambda \mbox{ for at most one } j=1,2,\dots \Biggr).
\end{eqnarray*}
Due to the assumption $1/2<\beta<1$, it follows that
\[
\varphi_n \Bigl( C_{\alpha}^{-1/\alpha} \epsilon\min
_{1 \leq i
\leq
m} \lambda_i \Bigr)\to0 %
\]
as $n\to\infty$; see \cite{samorodnitsky2004a}. Therefore, the proof
will be completed once we check that for all $\lambda>0$ and $0 \leq t
<t^\prime\leq1$,
\[
\psi_n\bigl(\lambda, t,t^\prime\bigr) \to0 %
\]
This, however, can be checked in exactly the same way as (4.10) in
\cite{owadasamorodnitsky2014}.
\end{pf}

\section{Other processes based on the range of the
subordinator} \label{secotherpr}

The distributional representation of the time-changed extremal process
(\ref{elimprocess}) in Corollary~\ref{cmainresult} can be stated
in the form
%
%
\begin{equation}
\label{ereprmain} Z_{\alpha,\beta}(t) = \hspace*{4pt}\rule{0pt}{13pt}^e\hspace*{-4pt}
\int_{(0,\infty)\times
\Omega^\prime}
\one \bigl( \bigl( R_\beta\bigl(\omega^\prime\bigr)+x \bigr)\cap
(0,t]\neq\varnothing \bigr) M\bigl(dx,d\omega^\prime\bigr),\qquad t\geq0.
\end{equation}
The self-similarity property of the process and the stationarity of
its max-increments can be traced to the scaling and shift invariance
properties of the range of the subordinator described in Proposition
\ref{prrange}. These properties can be used to construct other
self-similar processes with stationary max-increments, in the manner
similar to the way scaling and shift invariance properties of the real
line have been used to construct integral representations of Gaussian
and stable self-similar processes with stationary increments such as
fractional Brownian and stable motions; see,
for example, \cite{samorodnitskytaqqu1994} and
\cite{embrechtsmaejima2002}.

In this section, we describe one family of self-similar processes with
stationary max-increments, which can be viewed as an extension of the
process in~(\ref{ereprmain}). Other processes can be constructed; we
postpone a more general discussion to a later work.

For $0\leq s<t$, we define a function $j_{s,t}:  {\mathscr J} \to
[0,\infty]$ by
\[
j_{s,t}(F) = \sup \bigl\{ b-a: s<a<t, a,b\in F, (a,b)\cap F=\varnothing
\bigr\}, %
\]
the ``length of the longest empty space within $F$ beginning between
$s$ and $t$''. The function $j_{s,t}$
is continuous, hence measurable, on $\mathscr J$. Set also
$j_{s,s}(F)\equiv0$. Let
%
%
\begin{equation}
\label{egammar} 0<\gamma< (1-\beta)/\alpha,
\end{equation}
and define
%
%
\begin{eqnarray}\label{enewprocess}
Z_{\alpha,\beta,\gamma}(t) &=& \hspace*{4pt}\rule{0pt}{13pt}^e\hspace*{-4pt}
\int_{(0,\infty)\times
\Omega^\prime}
\bigl[\one \bigl( \bigl( R_\beta\bigl(\omega^\prime\bigr)+x \bigr)
\cap (0,t]\neq\varnothing \bigr) j _{0,t} \bigl(
R_\beta\bigl(\omega^\prime\bigr)+x \bigr) \bigr]^\gamma
\nonumber\\[-8pt]\\[-8pt]\nonumber
&&\hspace*{43pt}{}\times M \bigl(dx,d\omega^\prime\bigr),\qquad t\geq0.
\end{eqnarray}
It follows from~(\ref{eovershootd}) that
\[
E^\prime \biggl( \int_0^\infty \bigl[\one
\bigl( ( R_\beta+x )\cap (0,t]\neq\varnothing \bigr) j
_{0,t} ( R_\beta+x ) \bigr]^{\gamma\alpha} \beta x^{\beta-1}
\,dx \biggr)<\infty %
\]
for $\gamma$ satisfying~(\ref{egammar}). Therefore,
(\ref{enewprocess}) presents a well-defined Fr\'echet process. We
claim that this process is $H$-self-similar with
\[
H=\gamma+\beta/\alpha %
\]
and has stationary max-increments.

To check stationarity of max-increments, let $r>0$ and define
\begin{eqnarray}
Z^{(r)}_{\alpha,\beta,\gamma}(t) = \hspace*{4pt}\rule{0pt}{13pt}^e\hspace*{-4pt}
\int_{(0,\infty)\times
\Omega^\prime}
\bigl[\one \bigl( \bigl( R_\beta\bigl(\omega^\prime\bigr)+x \bigr)
\cap (r,r+t]\neq\varnothing \bigr) j _{r,r+t} \bigl(
R_\beta\bigl(\omega^\prime\bigr)+x \bigr)\bigr]^\gamma M
\bigl(dx,d\omega^\prime\bigr),\nonumber
\\
\eqntext{t\geq0.}
\end{eqnarray}
Trivially, for every $t\geq0$ we have
\[
Z_{\alpha,\beta,\gamma}(r) \vee Z^{(r)}_{\alpha,\beta,\gamma}(t) =
Z_{\alpha,\beta,\gamma}(r+t) %
\]
with probability 1, and it follows from part (c) of Proposition
\ref{prrange} that
\[
\bigl( Z^{(r)}_{\alpha,\beta,\gamma}(t), t\geq0 \bigr) \stackrel{d} {=}
\bigl( Z_{\alpha,\beta,\gamma}(t), t\geq0 \bigr). %
\]
Hence, stationarity of max-increments. Finally, we check the property
of self-similarity. Let $t_j>0,   \lambda_j>0,  j=1,\ldots, m$. Then
\[
P \bigl( Z_{\alpha,\beta,\gamma}(t_j) \leq\lambda_j, j=1,
\ldots, m \bigr)= \exp \bigl\{ -I(t_1,\ldots, t_m;
\lambda_1,\ldots, \lambda_m) \bigr\}, %
\]
where
\begin{eqnarray*}
&& I(t_1,\ldots, t_m; \lambda_1,\ldots,
\lambda_m) %
\\
&&\quad = E^\prime \biggl( \int_0^\infty\beta
x^{\beta-1} \max_{k=1,\ldots, m} \lambda_k^{-\alpha}
\bigl[\one \bigl( \bigl( R_\beta\bigl(\omega^\prime\bigr)+x
\bigr)\cap (0,t_k]\neq\varnothing \bigr) j _{0,t_k}
\bigl( R_\beta\bigl(\omega^\prime\bigr)+x \bigr) \bigr]^{\gamma\alpha}
\,dx \biggr). %
\end{eqnarray*}
Therefore, the property of self-similarity will follow once we check
that for any $c>0$,
\[
I(ct_1,\ldots, ct_m; \lambda_1,\ldots,
\lambda_m) = I\bigl(t_1,\ldots, t_m;
c^{-H}\lambda_1,\ldots, c^{-H}
\lambda_m\bigr). %
\]
This is, however immediate, since by using first part (b) of
Proposition~\ref{prrange} and, next, changing the variable of
integration to $y=x/c$ we have
\begin{eqnarray*}
&& I(ct_1,\ldots, ct_m; \lambda_1,\ldots,
\lambda_m)
\\
&&\quad = E^\prime \biggl(\int_0^\infty\beta
x^{\beta-1} \max_{k=1,\ldots,
m} \bigl\{ \lambda_k^{-\alpha}
\bigl[\one \bigl( \bigl( cR_\beta\bigl(\omega^\prime\bigr)+x
\bigr)\cap (0,ct_k]\neq\varnothing \bigr)
\\
&&\qquad {}\times \sup \bigl\{ b-a:
0<a<ct_j, a,b\in cR_\beta\bigl(\omega^\prime
\bigr)+x, (a,b)\cap cR_\beta\bigl(\omega^\prime\bigr)+x=\varnothing
\bigr\} \bigr]^{\alpha\gamma
} \bigr\} \,dx \biggr)
\\
&&\quad = c^{\beta+\alpha\gamma} E^\prime \biggl(\int_0^\infty
\beta x^{\beta-1} \max_{k=1,\ldots, m} \bigl\{ \lambda_k^{-\alpha}
\bigl[\one \bigl( \bigl( R_\beta\bigl(\omega^\prime\bigr)+x
\bigr)\cap (0,t_k]\neq\varnothing \bigr)
\\
&&\qquad{} \times  \sup \bigl\{ b-a:
0<a<t_k, a,b\in R_\beta\bigl(\omega^\prime\bigr)+x,
(a,b)\cap R_\beta\bigl(\omega^\prime\bigr)+x=\varnothing \bigr\}
\bigr]^{\alpha\gamma} \bigr\} \,dx \biggr)
\\
&&\quad = I\bigl(t_1,\ldots, t_m; c^{-H}
\lambda_1,\ldots, c^{-H}\lambda_m\bigr),
\end{eqnarray*}
as required.

\section*{Acknowledgments}
The anonymous referees read the paper very carefully, and some of
their comments were uncommonly perceptive and useful.

Samorodnitsky's research was partially supported by the ARO
Grant W911NF-12-10385 and NSA Grant H98230-11-1-0154 at Cornell
University and Fondation Mines Nancy.


%

\printhistory
\end{document}